\documentclass[12pt]{article}%
\usepackage[applemac]{inputenc}
\usepackage{amsmath,amssymb,fullpage}
\usepackage{amsfonts}
\usepackage{amsmath,amssymb}
\usepackage[applemac]{inputenc}
\usepackage{amsmath,amssymb,fullpage}
\usepackage{color}
\usepackage{color, colortbl}
\usepackage{amsmath}
\usepackage{amssymb}
\usepackage{graphicx}
\usepackage{tikz}
\usepackage{geometry}
\usetikzlibrary{arrows}
\usepackage{amsfonts}
\usepackage{amsmath,amssymb}
\usepackage[applemac]{inputenc}
\usepackage{amsmath,amssymb,fullpage}
\usepackage{color}
\usepackage{amsmath}
\usepackage{amssymb}
\usepackage{graphicx}
\usepackage{tikz}
\newtheorem{theorem}{Theorem}[section]

\newtheorem{proposition}[theorem]{Proposition}
\newtheorem{corollary}[theorem]{Corollary}
\newtheorem{definition}[theorem]{Definition}

\newtheorem{example}[theorem]{Example}

\newtheorem{problem}[theorem]{Problem}
\newtheorem{remark}[theorem]{Remark}

\usetikzlibrary{arrows,shapes,automata,petri}
\newenvironment{myenumerate}{

\begin{enumerate}}{\end{enumerate}}
\newcommand{\dproof}{\noindent {Proof.} \quad}
\newcommand{\fproof}{\hfill $\square$ \bigskip}

\numberwithin{equation}{section}

\definecolor{LightCyan}{rgb}{0.88,1,1}

\def\1B{\text{1\!\!I}}

\begin{document}

\title{Pricing of European options in incomplete jump diffusion  markets}
\author{Nacira Agram$^{1}$ \& Bernt \O ksendal$^{2}$}
\date{27 April 2021 \vskip 1.0cm
Dedicated to the memory of Mark H. A. Davis,\\
whose pioneering and innovative work inspired us all}
\maketitle

\footnotetext[1]{%
Department of Mathematics, Linnaeus University, V\" axj\" o, Sweden.
Email: nacira.agram@lnu.se. Work supported by the Swedish Research Council grant (2020-04697).}

\footnotetext[2]{%
Department of Mathematics, University of Oslo, Norway. 
Email: oksendal@math.uio.no.}
\begin{abstract}
We study option prices in financial markets where the risky asset prices are modelled by jump diffusions. Such markets are typically incomplete, and therefore there are in general infinitely many arbitrage-free option prices in these markets. We consider in particular European options with a terminal payoff $F$ at the terminal time $T$. It was proposed by Schweizer (1996) in a general semimartingale setting, following earlier works by F\" ollmer and Sondermann (1986) and Bouleau and Lamberton (1989), that the right price of such an option is the initial wealth needed to make it possible to generate by a self-financing portfolio a terminal wealth which is as close as possible to the payoff $F$ in the sense of variance. Schweizer calls this price \emph{the approximation price} and he investigates interesting general properties of this price and its corresponding optimal portfolio. 

However, neither of these authors compute explicitly this price in concrete cases. This is the motivation for the current paper: We apply stochastic control methods to compute this price in the setting of markets with assets described by jump diffusions.  Our method involves Stackelberg games and a suitably modified stochastic maximum principle. \\

We show that such an optimal initial wealth, denoted by $\widehat{z}$, with corresponding optimal portfolio $\widehat{\pi}$ exist and are unique.
We call $\widehat{z}$ the \emph{minimal variance price of $F$} and denote it by $p_{mv}(F)$. 

If the coefficients of the risky asset prices are deterministic, we show that 
$$p_{mv}(F)=E_{Q^{*}}[F],$$
for a specific equivalent martingale measure (EMM) $Q^{*}$. This shows in particular that the minimal variance price is free from arbitrage.

Then for the general case we apply a suitable maximum principle of optimal stochastic control to relate the minimal variance price  $\widehat{z}=p_{mv}(F)$ to the Hamiltonian and its adjoint processes, and we show that, under some conditions, $\widehat{z}=p_{mv}(F)=E_{Q_0}[F]$ for \emph{any} $Q_0$ in a family $\mathbb{M}_0$ of EMMs, described by the set of solutions of a system of linear equations.

Finally, we illustrate our results by looking at specific examples.

\end{abstract}



\vskip 0.5cm


\textbf{2010 AMS subject classifications:} Primary: 60G51; 60H10; 60J75; 93E20; 91Gxx; 60H35. Secondary: 60H30.
\newline

\textbf{Keywords :} Option pricing; Arbitrage-free; Equivalent martingale measure; Jump diffusion; Stochastic control.



\section{Introduction}
\quad The recent economic crises, including impacts from the Covid 19, has given researchers in mathematical finance a motivation to study new models for financial markets. In particular, it has become more relevant than ever to include default time or jumps that come from L\' evy  processes. However, such markets are typically incomplete and there are infinitely many EMMs and hence also infinitely many arbitrage-free prices for European options. Therefore it is natural to ask what is the best arbitrage-free option pricing rule to use in such incomplete markets? Or, equivalently, if $F$ is the terminal payoff of the option (and the interest rate is 0), what is the best EMM $Q$ to choose for the  pricing rule that
$$\text{"price of option"}=E_Q[F]?$$
This is a question that has been studied by meany researchers, and several approaches and answers have been provided. For example, we mention the following:

F\"ollmer and Schweizer (1991)  
discuss optimal hedging in incomplete markets and introduce the concept of a risk-minimising martingale measure, which they relate to optimality of strategies, but not directly to option pricing.

Delbaen and Schachermayer (1996)  
consider continuous markets and propose to use as the option pricing measure the EMM obtained by minimising the $L^2$-norm of its Radon-Nikodym derivative with respect to the original probability measure $P$. 

Frittelli (2000)  
deals with a general market model and suggests to use the EMM with minimal entropy. He justifies this by relating it to the marginal utility of terminal wealth in an exponential utility maximisation problem. Sufficient conditions for a martingale measure to be entropy-minimising are given by Grandits and Rheinlaender (2002) , using the theory of BMO-martingales.

Hodges and Neuberger (1989) 
introduce the utility indifference pricing principle, saying that the price $p$ at time 0 of an option with terminal payoff $F$ should be such that a seller of the option is indifferent, with respect to a given utility function, to the following two options:\\
(i) either receiving that price $p$ for the option and at time 0, then trade optimally with that added initial wealth and paying out the amount $F$ at the terminal time $T$, or\\
(ii) not selling the option at all and just trade optimally without any payout at the terminal time $T$.
There are many papers dealing with the use of this pricing principle. For a survey see Henderson and Hobson (2004)  
and the references therein. In the case of exponential utility, this pricing rule turns out to be related to pricing by means of the minimal entropy measure. See Davis and Yoshikawa (2016). 

It was proposed by Schweizer (1996), in a general semimartingale setting, following earlier works by F\" ollmer and Sondermann (1986) and Bouleau and Lamberton (1989), that the right price of such an option is the initial wealth needed to make it possible to generate by a self-financing portfolio a terminal wealth which is as close as possible to the payoff $F$ in the sense of variance. 
This pricing principle is natural in view of the well-known result that if the terminal payoff $F$ of the European option is replicable, then the unique arbitrage-free option price of the option is the initial wealth $z=X(0)$ needed for the replication, i.e. to make $X(T)=F$ a.s., where $X(t)$ denotes the wealth at time $t$, $0 \leq t \leq T$. Therefore one can argue that in general the right price of a European option with terminal payoff $F$ (replicable or not) should be the initial wealth $z=X(0)$ needed to make it possible to generate by means of an admissible portfolio $\pi$ a terminal wealth $X_{z,\pi}(T)$ which is as close as possible to $F$ in the sense of variance,
i.e. $\widehat{z},\widehat{\pi}$ minimizes 
$$(z,\pi) \mapsto E[(X_{z,\pi}(T)-F)^2],$$
where $E[\cdot]=E_{P}[\cdot]$ denotes expectation with respect to the underlying probability measure $P$.

Schweizer calls this price \emph{the approximation price} and he investigates interesting general properties of this price and its corresponding optimal portfolio. 

However, neither of these authors compute explicitly this price in concrete cases. This is the motivation for the current paper: We apply stochastic control methods to compute this price in the setting of markets with assets described by jump diffusion.  Our method involves Stackelberg games and a suitably modified stochastic maximum principle. \\
Our paper is organised as follows:

In Section 2 we show that such an optimal initial wealth $z=\widehat{z}$ with corresponding optimal portfolio $\pi=\widehat{\pi}$ exist and are unique.
We call $\widehat{z}$ the \emph{minimal variance price of $F$} and denote it by $p_{mv}(F)$. In the classical Black-Scholes market this price coincides with the classical Black-Scholes option price. To find this optimal initial wealth $\widehat{z}$  we use two methods:

 In the first method, in Section 3 we formulate the problem as a Stackelberg game and find its unique solution using stochastic calculus. In particular, we show that under some conditions there is a unique EMM $Q^{*}$ which gives the solution $\widehat{z}$ of the minimal variance problem, in the sense that 
$\widehat{z}=p_{mv}(F)=E_{Q^{*}}[F]$. This shows in particular that the minimal variance price is free from arbitrage.

In the second method, in Section 4 we introduce a suitable version of the maximum principle for optimal stochastic control to the problem. Using this we can relate the EMMs $Q$ such that $\widehat{z}=E_{Q}[F]$ to the Hamiltonian and the corresponding adjoint processes $p(t),q(t),r(t,\zeta)$. Then we show that, under some conditions, $\widehat{z}=p_{mv}(F)=E_{Q_0}[F]$ for \emph{any} $Q_0$ in a family $\mathbb{M}_0$ of EMMs, described by the set of solutions of a system of linear equations.

Finally, in Section 5 we give examples to illustrate our results.

\section{Option pricing in the general incomplete market case}
In this section we present our financial market model, and we give a brief survey of some of the fundamental concepts and results from the theory of pricing of European options. 

Consider a financial market with two investment possibilities:

\begin{description}
\item[(i)] A risk free asset, with unit price $S_{0}(t)=1$ for all $t$,

\item[(ii)] A risky asset, with unit price $S(t)$ at time $t$ given by, writing $\mathbb{R}^{*}=\mathbb{R}\setminus \{0\}$,%
\begin{align}\label{S}
dS(t) & =S(t^{-})\left[  \alpha(t)dt+\sigma(t)dB(t)+\int_{
\mathbb{R^*}
}\gamma(t,\zeta)\widetilde{N}(dt,d\zeta)\right],\quad S(0) >0.
\end{align}
\end{description}
For simplicity of notation, we will write $S(t)$ instead of $S(t^{-})$ in the following.\\
Here $\alpha(t) \in \mathbb{R}, \sigma(t)=(\alpha_1(t), ..., \alpha_m(t)) \in \mathbb{R}^{m},$ $\gamma(t,\zeta)=(b_1(t,\zeta),..., b_k(t,\zeta)) \in \mathbb{R}^{k}$, and $B(t)=(B_1(t), ... , B_m(t))^{'} \in \mathbb{R}^{m}$ and  
$\widetilde{N}(dt,d\zeta)=(\widetilde{N}_1(dt,d\zeta), ..., \widetilde{N}_k(dt,d\zeta))^{'} \in \mathbb{R}^{k}$ are independent Brownian motions and compensated Poisson random measures, respectively, on a complete filtered probability space $(\Omega,\mathcal{F},\mathbb{F}=\{\mathcal{F}_t\}_{t\geq 0},P)$. We are using the matrix notation, i.e.
\begin{align*}
 \sigma(t)dB(t):&=\sum_{i=1}^{m} \sigma_i(t)dB_i(t),\text{  }
 \gamma(t,\zeta)\widetilde{N}(dt,d\zeta)= \sum_{j=1}^{k}\gamma_j(t,\zeta)\widetilde{N}_j(dt,d\zeta),
\end{align*}
and we assume that $\gamma(t,\zeta)>-1$ and 
\begin{equation}\label{1.2}
\sigma^2(t) +\int_{
\mathbb{R^*}
} \gamma^2(t,\zeta) \nu(d\zeta) > 0 \text{ for a.a. } t, \text{ a.s.}
\end{equation}
For simplicity we assume throughout this paper that all the coefficients $\alpha,\sigma$ and $\gamma$ are bounded $\mathbb{F}$-predictable processes.
The results in this paper can easily be extended to an arbitrary number of risky assets, but since the features of incomplete markets we are dealing with, can be fully illustrated by one jump diffusion risky asset only, we will for simplicity concentrate on this case in the following. We emphasise however, that we deal with an arbitrary number $m$ of independent Brownian motions and an arbitrary number $k$ of independent Poisson random measures in the representation \eqref{S}.

Let $z\in%
\mathbb{R}
$ be a given initial wealth and let $\pi(t)\in \mathbb{R}$ be a self-financing portfolio,
representing the fraction of the total wealth $X(t)=X_{z,\pi}(t)$ invested in
the risky asset at time $t$. Then the corresponding wealth dynamics is given
by the following linear stochastic differential equation (SDE) with jumps
\begin{align}\label{W}
dX(t) =X(t)\pi(t)\left[  \alpha(t)dt+\sigma(t)dB(t)+\int_{
\mathbb{R^*}
}\gamma(t,\zeta)\widetilde{N}(dt,d\zeta)\right], \quad X(0) =z.
\end{align}
Let $\mathbb{F}=\{\mathcal{F}_{t}\}_{t\geq0}$ be the filtration generated by
$\{B(s)\}_{s\leq t}$ and $\{N(s,\zeta)\}_{s\leq t}$. \vskip 0.2cm

For the convenience of the reader, let us recall the general pricing problem (of a European option):\\
Let $F$ be a given \emph{T-claim}, i.e. $F \in L^2(P)$ is an  $\mathcal{F}_{T}$-measurable random variable, representing the payoff at the terminal time $T$ written on the contract. 
If you own the contract, you are entitled to get this payoff $F$ at time $T$. For example, $F$ might be linked to the terminal value of the risky asset, i.e. $F=h(S(T))$ for some function $h$. If somebody comes to you and offers you this contract, how much would you be willing to pay for it now, at time $t=0$? \\
To answer this question the \emph{buyer}, who does not want to risk losing money, will argue as follows:\\
If the price of the option is $z$, then to buy it I need to borrow $z$ in the bank, thereby starting with an initial wealth $-z$. Then it should be possible for me to find a portfolio $\varphi$ such that the corresponding final wealth $X_{-z,\varphi}(T)$ with the contract payoff $F$ added, gives me a non-negative net wealth, a.s. In other words, the buyer's price of $F$, denoted by $p_{b}(F)$, is defined by
\[
p_{b}(F)=\sup\{z\text{; there exists }\varphi\in\mathcal{A}\text{ such that
}X_{-z,\varphi}(T)+F\geq0\text{ a.s.}\}.
\]
Similarly, the \emph{seller} of the contract does not want the risk of losing money either, so he/she argues as follows:
If I receive a payment $y$ for the contract, it should be possible for me to find a portfolio $\psi $, such that the corresponding wealth at time $T$ is big enough to pay the contract obligation $F$ at the terminal time $T$, a.s.
\\
In other words, the seller's price of $F$, denoted by $p_{s}(F)$, is defined by
\[
p_{s}(F)=\inf\{y:\text{ there exists }\psi\in\mathcal{A}\text{ such that
}X_{y,\psi}(T)\geq F\text{ a.s.}\}.
\]
In general the two prices $p_b(F)$ and $p_{s}(F)$ do not coincide. In fact, we have the following well-known result:

Let $\mathbb{M}$ denote the set of all EMMs for $S(\cdot)$.\\
Then%
\begin{equation}
p_{b}(F)\leq E_{Q}[F]\leq p_{s}(F), \label{pz}%
\end{equation}
for all $Q\in\mathbb{M}$. 
See e.g. 
\O ksendal \& Sulem (2019) \cite{OS}.

A portfolio $\pi$ is called an \emph{arbitrage}  if it can generate from an initial wealth $X(0)=z=0$ a terminal wealth $X_{0,\pi}(T) \geq 0$ a.s. with $P[X_{0,\pi}(T) > 0] > 0.$ In that sense, an arbitrage is a kind of money machine; it can generate a profit with no risk.  Markets with arbitrage cannot exist for a reasonable length of time. Therefore option prices which lead to arbitrage are not accepted in any financial market model. An option  price $p(F) \in \mathbb{R}$ is called \emph{arbitrage-free} if it does not give an arbitrage to neither the buyer nor the seller. It is well-known and easy to prove that for any $Q \in \mathbb{M}$ the price $E_Q[F]$ is arbitrage-free.\\

Even though the buyer's price and the seller's price are different in general, they can be the same in some cases. For example, it is well-known that if $F$ is \emph{replicable}, in the sense that there exist $z\in\mathbb{R}$ and $\pi\in\mathcal{A}$ (the set of admissible portfolios), such that
\begin{align*}
X_{z,\pi}(T)=F\text{ a.s.,}%
\end{align*}
then we have \emph{equality} in (\ref{pz}), i.e,%
\begin{align}\label{1.5}
p_{b}(F)=E_{Q}[F]= p_{s}(F),\text{ for \bf{all} } Q\in\mathbb{M}.
\end{align}
This is a special case, because in a general incomplete market there are many payoffs  $F$ which are not replicable. Only if the market is complete, are all $F$  replicable, and in that case there is only one EMM $Q$.
In the next section we will prove a generalisation of \eqref{1.5}. See Theorem \ref{th2.2} and Corollary \ref{cor2.3}.

Thus the interesting case is when the market is incomplete, and  then there is usually a big gap between $p_s$ and $p_b$. Then one needs additional requirements to find the "best" price.
The purpose of this paper is to use a \emph{stochastic control approach} to
compute the approimate/minimal variance price of Schweizer explicitly for markets described by jump diffusion processes. Specifically, if $F$ is a given T-claim, then to each initial wealth $z\in%
\mathbb{R}
$ and each portfolio $\pi \in\mathcal{A}$, we associate the quadratic
cost functional%
\begin{align}\label{cost}
 J(z,\pi)=E\left[  \frac{1}{2}\left(  X_{z,\pi}(T)-F\right)  ^{2}\right]  .   
\end{align}
Then we consider the following problem:

\begin{problem} \label{mv}
Find the optimal initial wealth $\widehat{z}\in%
\mathbb{R}
$ and the optimal portfolio $\widehat{\pi}\in\mathcal{A}$,  such that%
\begin{equation}
\underset{z,\pi}{\inf}J(z,\pi)=J(\widehat{z},\widehat{\pi}). \label{cost}%
\end{equation}
\end{problem}
\emph{Heuristically, this means that we define the price of the option with payoff $F$ to be the initial wealth $\widehat{z}$ needed to get the terminal wealth $X(T)$ as close as possible to $F$ in quadratic mean by an admissible portfolio.} 

In other words, we make the following definition:

\begin{definition}
Suppose that $(\widehat{z},\widehat{\pi})\in%
\mathbb{R}
\times\mathcal{A}$ is the unique solution of the stochastic control problem Problem \ref{mv}.
Then we define \emph{the minimal variance price} of $F$, denoted by $p_{mv}(F)$,  by
\begin{equation*}
p_{mv}(F)=\widehat{z}.
\end{equation*}
\end{definition}

\section{Existence and uniqueness of the optimal initial wealth}
In this section we prove the existence and the uniqueness of solutions of Problem \ref{mv}. 
\subsection{Equivalent martingale measures (EMMs)}
Since an EMMs play a crucial role in our discussion, we start this section by recalling that an important group of measures $Q\in\mathbb{M}$ can be described as follows (we refer to Chapter 1 in \O ksendal \& Sulem (2019) 
 for more details):\\
Let $\theta_{0}(t)$ and $\theta_{1}(t,\zeta)>-1$ be $\mathbb{F}$-predictable
processes such that
\begin{align}\label{emm}
\alpha(t)+\theta_{0}(t)\sigma(t)+\int_{
\mathbb{R^*}
}\theta_{1}(t,\zeta)\gamma(t,\zeta)\nu(d\zeta)=0\text{, } \quad t\geq0.
\end{align}
Define the local martingale $Z(t)=Z^{\theta_{0},\theta_{1}}(t)$, by%

\begin{align}\label{Z1}
dZ(t)  =Z(t)\left[  \theta_{0}(t)dB(t)+\int_{\mathbb{R^*}
}\theta_{1}(t,\zeta)\widetilde{N}(dt,d\zeta)\right], \quad 
Z(0) =1,
\end{align}
i.e.,
\small
\begin{align}\label{Z2}
Z(t)  &  =\exp\left(  \int_{0}^{t}\theta_{0}(s)dB(s)-\frac{1}{2}\int_{0}%
^{t}\theta_{0}^{2}(s)ds\right.  +\int_{0}^{t}\int_{\mathbb{R^*}
}\left\{  \ln(1+\theta_{1}(s,\zeta))-\theta_{1}(s,\zeta)\right\}  \nu
(d\zeta)ds\nonumber\\
&  \left.  +\int_{0}^{t}\int_{\mathbb{R^*}
}\ln(1+\theta_{1}(s,\zeta))\widetilde{N}(ds,d\zeta)\right)  .
\end{align}
Suppose that $Z$ is a true martingale.  A sufficient condition for this to hold is
\begin{align}\label{KS}
    E\Big[\exp\Big(\frac{1}{2}\int_0^T \theta_0^2(s)ds + \int_0^T\int_{\mathbb{R^*}}\theta_1^2(s,\zeta) N(ds,d\zeta)\Big)\Big] < \infty.
\end{align}
See Kallsen and Shiryaev (2002).
Then the measure $Q^{\theta_{0},\theta_{1}}$ defined by%
\begin{align} \label{Q}
dQ^{\theta_{0},\theta_{1}}(\omega)=Z(T)dP(\omega)\text{ on }\mathcal{F}_{T}%
\end{align}
is in $\mathbb{M}$.

\subsection{The optimal portfolio}
We may regard the minimal variance problem (Problem \ref{mv}) as a Stackelberg game, in which the first player chooses the initial wealth $z$, followed by the second player choosing the optimal portfolio $\pi$ based on this initial wealth. Knowing this response $\pi=\pi_{z}$ from the follower, the first player chooses the initial wealth $\widehat{z}$ which leads to a response $\pi=\pi_{\widehat{z}}$ which is optimal, in the sense that $J(\widehat{z},\pi_{\widehat{z}}) \leq J(z,\pi)$ over all admissible pairs $(z,\pi)$.
To this end, in this subsection we first proceed to find the optimal portfolio based on a given initial wealth $z$. \\
Accordingly, assume as before that the wealth process $X(t)=X_{z, \pi}(t)$, corresponding to an initial wealth $z$ and a self-financing portfolio $\pi$, is given by
\begin{align}\label{3.1a}
    dX(t)=X(t)\pi(t)\Big[\alpha(t) dt + \sigma(t) dB(t) + \int_{ \mathbb{R^*} }\gamma(t,\zeta)\widetilde{N}(dt,d\zeta)\Big],\quad
    X(0)=z.
\end{align}
Suppose the terminal payoff $F\in L^2(P,\mathcal{F}_T)$ has the form $F=F(T)$, where the martingale $F(t):=E[F|\mathcal{F}_t], t \in [0,T]$ has the  It\^o-L\'evy representation
\begin{align*}
    dF(t)=\beta(t)dB(t) + \int_{ \mathbb{R^*} } \kappa(t,\zeta) \widetilde{N}(dt,d\zeta), \quad E[F]= F_0.
\end{align*}
for some (unique) $\mathbb{F}$-predictable processes $\beta(t) \in L^2(\lambda \times P),
\kappa(t,\zeta) \in L^2(\lambda \times \nu \times P)$.
Then by the It\^{o} formula for jump diffusions (see e.g Theorem 1.14 in \O ksendal and Sulem (2019))  \cite{OS}
we get
\begin{align}
    d(X(t)F(t))&= X(t) dF(t) + F(t) dX(t) + d[X,F]_t \nonumber\\
    &= X(t)[\beta(t)dB(t) + \int_{ \mathbb{R^*} } \kappa(t,\zeta)\widetilde{N}(dt,d\zeta)]\nonumber\\
    &+X(t)F(t)[\pi(t)\alpha(t) dt + \pi(t) \sigma(t) dB(t) + \int_{ \mathbb{R^*} } \pi(t)\gamma(t,\zeta)\widetilde{N}(dt,d\zeta)]\nonumber\\
    & + X(t)[\pi(t)\sigma(t)\beta(t) dt + \int_{ \mathbb{R^*} } \pi(t)\gamma(t,\zeta)\kappa(t,\zeta)\widetilde{N}(dt,d\zeta)\nonumber\\
 &+  \int_{ \mathbb{R^*} } \pi(t)\gamma(t,\zeta)\kappa(t,\zeta)\nu(d\zeta)dt.] \label{3.3a}
\end{align}
Hence
\begin{align}
    E[X(T)F(T)]&=z F_0+\int_0^T E\Big[ X(t)\Big\{F(t)\pi(t)\alpha(t)
    + \pi(t)\sigma(t)\beta(t) \nonumber \\
    &+ \int_{ \mathbb{R^*} } \pi(t)\gamma(t,\zeta)\kappa(t,\zeta)\nu(d\zeta)\Big\}\Big]dt. \label{3.4a}
\end{align}
Similarly,
\begin{align*}
E[X^2(T)]= z^2 + \int_0^T E\Big[ X^2(t)\Big\{2 \pi(t)\alpha(t) + 
\pi^2(t)\sigma^2(t) + \int_{ \mathbb{R^*} } \pi^2(t) \gamma^2(t,\zeta) \nu(d\zeta)\Big\} \Big]dt,
\end{align*}
and
\begin{align*}
    E[F^2(T)]=F_0^2+\int_0^TE\Big[\beta^2(t) + \int_{ \mathbb{R^*} }\kappa^2(t,\zeta) \nu(d\zeta)\Big]dt. 
\end{align*}
This gives
\small
\begin{align*}
    J(\pi)&=E\Big[\frac{1}{2} (X(T) -F)^2\Big]= \frac{1}{2}\Big(E[X^2(T)] - 2E[X(T)F(T)]+E[F^2(T)]\Big)\nonumber\\
    &=\frac{1}{2}(z-F_0)^2 + \frac{1}{2}E\Big[\int_0^T\Big\{2 X^2(t)\pi(t)\alpha(t)
    +  X^2(t)\pi^2(t) \Big(\sigma^2(t)+\int_{ \mathbb{R^*} } \gamma^2(t,\zeta) \nu(d\zeta)\Big)\nonumber\\
    &-2X(t)\pi(t)\Big(F(t) \alpha(t)
    + \sigma(t)\beta(t) +\int_{ \mathbb{R^*} }\gamma(t,\zeta)\kappa(t,\zeta) \nu(d\zeta)\Big)\nonumber\\ 
    &+\beta^2(t) + \int_{ \mathbb{R^*} }\kappa^2(t,\zeta) \nu(d\zeta)\Big\} dt\Big].
\end{align*}
We can minimise $J(\pi)$ by minimising the $dt$-integrand pointwise for each $t$. This gives the following result:
\begin{theorem} \label{th2.1}
a) For given initial value $X(0)=z > 0$ the portfolio $\widehat{\pi}=\widehat{\pi}_{z}$ which minimises
\begin{align*}
   \pi \mapsto  E\Big[\frac{1}{2}(X_{z,\pi}(T)-F)^2\Big]
\end{align*}
is given in feedback form with respect to $X(t)=X_{z,\widehat{\pi}}$ by
\begin{align}
    \widehat{\pi}(t)&=\widehat{\pi}(t,X(t))= \frac{F(t)\alpha(t) +\sigma(t)\beta(t)+\int_{ \mathbb{R^*} }\gamma(t,\zeta)\kappa(t,\zeta)\nu(d\zeta) - X(t)\alpha(t)}
    {X(t) \Big( \sigma^2(t) + \int_{ \mathbb{R^*} } \gamma^2(t,\zeta) \nu(d\zeta) \Big)}, \label{3.9a}
\end{align}
or, equivalently,
\begin{align}\label{pi}
    \widehat{\pi}(t)X(t)=G(t)[X(t)-F(t)]+\frac{\sigma(t)\beta(t)+\int_{\mathbb{R^*}}\gamma(t,\zeta)\kappa(t,\zeta)\nu(d\zeta)}{\sigma^2(t)+\int_{\mathbb{R^*}}\gamma^2(t,\zeta)\nu(d\zeta)},
\end{align}
where
\begin{equation} \label{(6)}
    G(t)=-\alpha(t)\Big(\sigma^2(t)+\int_{\mathbb{R^*}} \gamma^2(t,\zeta) \nu(d\zeta)\Big)^{-1}.
\end{equation}
b) Given an initial value $z>0$, the corresponding optimal wealth $X_{\widehat{\pi}}(t)=\widehat{X}(t)$ solves the SDE 
\begin{align} \label{2.10}
    d\widehat{X}(t)&=\widehat{X}(t)\widehat{\pi}(t,\widehat{X}(t))\Big[\alpha(t) dt + \sigma(t) dB(t) + \int_{ \mathbb{R^*} }\gamma(t,\zeta)\widetilde{N}(dt,d\zeta)\Big]\nonumber\\
   &=\frac{F(t)\alpha(t) +\sigma(t)\beta(t)+\int_{ \mathbb{R^*} }\gamma(t,\zeta)\kappa(t,\zeta)\nu(d\zeta) - X(t)\alpha(t)}
    {\sigma^2(t) + \int_{ \mathbb{R^*} } \gamma^2(t,\zeta) \nu(d\zeta)}\times\nonumber\\
    &\times \Big[\alpha(t) dt + \sigma(t) dB(t)+\int_{ \mathbb{R^*} }\gamma(t,\zeta)\widetilde{N}(dt,d\zeta)\Big].
\end{align}
\end{theorem}
\begin{remark}
\textbf{(a)} Consider the special case when $F$ is a deterministic constant. Then $F(t)=F=E[F]$ for all $t$, and $\beta=\kappa =0$. Hence the optimal portfolio is given in feedback form by
$$\widehat{\pi}(t)X(t)=G(t)[X(t)-F].$$
Assume, for example, that $\alpha(t) > 0$. Then $G(t) < 0$ and we see that if $X(t) < F$ then $\widehat{\pi}X(t) >0$ and hence the optimal portfolio pushes $X(t)$ upwards towards F. Similarly, if $X(t) > F$ then $\widehat{\pi}(t)X(t) < 0$ and the optimal push of $X(t)$  is downwards towards $F$. This is to be expected, since the portfolio tries to minimise the terminal variance $E[(X(T)-F)^2]$.
\vskip 0.2cm
\textbf{(b)} In particular, we see that if $F$ is a deterministic constant, and we start at $z=X(0)=F$, we can choose $\pi=0$ and this gives 
$J(z,\pi)=J(F,0)=E[\frac{1}{2}(X(T)-F)^2]=E[\frac{1}{2}(F-F)^2]=0$, which is clearly optimal. By uniqueness of $(\widehat{z},\widehat{\pi})$ we conclude that $(\widehat{z},\widehat{\pi})=(F,0)$ is the optimal pair in this case.  
\end{remark}

\subsection{Explicit expression of the corresponding optimal wealth 
}

 Writing $X=\widehat{X}$ for notational simplicity, equation \eqref{2.10} is of the form
\begin{equation}\label{(1)}
dX\left( t\right) = C\left( t\right) d\Lambda _{t}+X\left( t\right)
d\Gamma _{t} , 
X\left( 0\right) = z,
\end{equation}%
where%
\begin{equation}
C\left( t\right) =\frac{F\left( t\right) \alpha \left( t\right) +\sigma
\left( t\right) \beta \left( t\right) +\int_{ \mathbb{R^*} }\gamma \left(
t,\zeta\right) \kappa\left( t,\zeta\right) \nu \left( d\zeta\right) }{ \sigma
^{2}\left( t\right) +\int_{ \mathbb{R^*} }\gamma ^{2}\left( t,\zeta\right) \nu
\left( d\zeta\right) },  \label{(2)}
\end{equation}%
\begin{equation}
d\Lambda _{t}=\alpha \left( t\right) dt+\sigma \left( t\right) dB\left(
t\right) +\int_{ \mathbb{R^*} }\gamma \left( t,\zeta\right) \widetilde{N}\left(
dt,d\zeta\right),  \label{(3)}
\end{equation}%
\begin{equation}
d\Gamma _{t}=\alpha _{1}\left( t\right) dt+\sigma _{1}\left( t\right)
dB\left( t\right) +\int_{ \mathbb{R^*} }\gamma _{1}\left( t,\zeta\right) \widetilde{N%
}\left( dt,d\zeta\right) ,  \label{(4)}
\end{equation}%
with%
\begin{equation}
\alpha _{1}\left( t\right) =G\left( t\right) \alpha \left( t\right) ,\sigma
_{1}\left( t\right) =G\left( t\right) \sigma \left( t\right) ,\gamma
_{1}\left( t,\zeta\right) =G\left( t\right) \gamma \left( t,\zeta\right).
\label{(5)}
\end{equation}%
We rewrite (\ref{(1)}) as%
\begin{equation}
dX\left( t\right) -X\left( t\right) d\Gamma _{t}=C\left( t\right) d\Lambda
_{t}, \label{(7)}
\end{equation}%
and multiply this equation by a process of the form%
\begin{equation}
Y_{t}=Y_{t}^{\left( \rho ,\lambda ,\theta \right) }=\exp \left(
A_{t}^{\left( \rho ,\lambda ,\theta \right) }\right),  \label{(8)}
\end{equation}%
with%
\begin{equation}
A_{t}^{\left( \rho ,\lambda ,\theta \right) }=\int_{0}^{t}\rho \left(
s\right) ds+\int_{0}^{t}\lambda \left( s\right) dB\left( s\right)
+\int_{0}^{t}\int_{ \mathbb{R^*} }\theta \left( s,\zeta\right) \widetilde{N}\left(
ds,d\zeta\right) ,  \label{(9)}
\end{equation}%
where $\rho ,\lambda $ and $\theta $ are processes to be determined.\\
Then (\ref{(7)}) gets the form%
\begin{equation}
Y_{t}dX\left( t\right) -Y_{t}X\left( t\right) d\Gamma _{t}=Y_{t}C\left(
t\right) d\Lambda _{t}.  \label{(10)}
\end{equation}%
We want to choose $\rho ,\lambda $ and $\theta $ such that $Y_t$ becomes an integrating factor, in the sense that
\begin{equation}
d\left( Y_{t}X\left( t\right) \right) =Y_{t}dX\left( t\right) -Y_{t}X\left(
t\right) d\Gamma _{t}+\text{terms not depending on }X. \label{(11)}
\end{equation}%
To this end, note that by the It\^{o} formula for L\' evy processes, we have
\begin{eqnarray}
dY_{t} &=&Y_{t}\left[ \rho \left( t\right) dt+\lambda \left( t\right)
dB\left( t\right) \right] +\frac{1}{2}Y_{t}\lambda ^{2}\left( t\right) dt
\label{(12)} \\
&&+\int_{ \mathbb{R^*} }\left\{ \exp \left( A_{t}+\theta \left( t,\zeta\right)
\right) -\exp \left( A_{t}\right) -\exp \left( A_{t}\right) \theta \left(
t,\zeta\right) \right\} \nu \left( d\zeta\right) dt  \nonumber \\
&&+\int_{ \mathbb{R^*} }\left\{ \exp \left( A_{t}+\theta \left( t,\zeta\right)
\right) -\exp \left( A_{t}\right) \right\} \widetilde{N}\left( dt,d\zeta\right) 
\nonumber \\
&=&Y_t \left[ \left\{ \rho \left( t\right) +\frac{1}{2}\lambda ^{2}\left(
t\right) +\int_{ \mathbb{R^*} }\left( e^{\theta \left( t,\zeta\right) }-1-\theta
\left( t,\zeta\right) \right) \nu \left( d\zeta\right) \right\} dt\right.   \nonumber
\\
&&\left. +\lambda \left( t\right) dB\left( t\right) +\int_{ \mathbb{R^*} }\left(
e^{\theta \left( t,\zeta\right) }-1\right) \widetilde{N}\left( dt,d\zeta\right)
\right] .  \nonumber
\end{eqnarray}%
Therefore, again by the It\^{o} formula, using \eqref{(1)},
\begin{eqnarray}
d\left( Y_{t}X\left( t\right) \right)  &=&Y_{t}dX\left( t\right) +X\left(
t\right) dY_{t}+d\left[ X,Y\right] _{t}  \label{(13)} \\
&=&Y_{t}dX\left( t\right) +Y_{t}X\left( t\right) \left[ \left\{ \rho +\frac{1%
}{2}\lambda ^{2}+\int_{ \mathbb{R^*} }\left( e^{\theta \left( t,\zeta\right)
}-1-\theta \left( t,\zeta\right) \right) \nu \left( d\zeta\right) \right\} dt\right. 
\nonumber \\
&&\left. +\lambda \left( t\right) dB\left( t\right) +\int_{ \mathbb{R^*} }\left(
e^{\theta \left( t,\zeta\right) }-1\right) \widetilde{N}\left( dt,d\zeta\right) %
\right]   \nonumber \\
&&+Y_{t}X\left( t\right) \left[ \left\{ \lambda \left( t\right) \sigma
_{1}\left( t\right) +\int_{ \mathbb{R^*} }\left( e^{\theta \left( t,\zeta\right)
}-1\right) \gamma _{1}\left( t,\zeta\right) \nu \left( d\zeta\right) \right\}
dt\right.   \nonumber \\
&&\left. +\int_{ \mathbb{R^*} }\left( e^{\theta \left( t,\zeta\right) }-1\right)
\gamma _{1}\left( t,\zeta\right) \widetilde{N}\left( dt,d\zeta\right) \right]
+Y_{t}C\left( t\right) dK_{t},  \nonumber
\end{eqnarray}%
where%
\begin{equation}
dK_{t}=\left\{ \lambda \left( t\right) \sigma \left( t\right) +\int_{ \mathbb{R^*} }\left( e^{\theta \left( t,\zeta\right) }-1\right) \gamma \left( t,\zeta\right)
\nu \left( d\zeta\right) \right\} dt+\int_{ \mathbb{R^*} }\left( e^{\theta \left(
t,\zeta\right) }-1\right) \gamma \left( t,\zeta\right) \widetilde{N}\left(
dt,d\zeta\right).  \label{(14)}
\end{equation}%
This gives 
\begin{align}
&d\left( Y_{t}X\left( t\right) \right) -Y_{t}dX\left( t\right)
+Y_{t}X\left( t\right) d\Gamma _{t}  \label{(15)} \\
&=Y_{t}X\left( t\right) \Big[ \Big\{ \rho +\alpha _{1}+\frac{1}{2}\lambda
^{2}+\lambda \sigma_1+\int_{ \mathbb{R^*} }\left( e^{\theta \left( t,\zeta\right) }-1-\theta \left(
t,\zeta\right) \right) \nu \left( d\zeta\right) \Big\}dt   \nonumber \\
&+ \left( \lambda \left( t\right) +\sigma _{1}\left( t\right) \right)
dB\left( t\right) +\int_{ \mathbb{R^*} }\left( e^{\theta \left( t,\zeta\right)
}-1\right) \gamma _{1}\left( t,\zeta\right) \nu \left( d\zeta\right) dt 
\nonumber \\
& +\int_{ \mathbb{R^*} }\Big\{ \left( e^{\theta \left( t,\zeta\right)
}-1\right) \left( 1+\gamma _{1}\left( t,\zeta\right) \right) +\gamma _{1}\left(
t,\zeta\right) \Big\} \widetilde{N}\left( dt,d\zeta\right)   \nonumber \\
&+Y_{t}C\left( t\right) 
dK_t. \nonumber
\end{align}%
Choose $\theta \left( t,\zeta\right) =\widehat{\theta }\left( t,\zeta\right) $, such
that 
\[
( e^{\theta \left( t,\zeta\right) }-1) \left( 1+\gamma _{1}\left(
t,\zeta\right) \right) +\gamma _{1}\left( t,\zeta\right) =0,
\]%
i.e.
\begin{equation}
\widehat{\theta }\left( t,\zeta\right) =-\ln \left( 1+\gamma _{1}\left(
t,\zeta\right) \right).  \label{(16)}
\end{equation}%
Next, choose $\lambda \left( t\right) =\widehat{\lambda }\left( t\right) $
such that 
\begin{equation}
\widehat{\lambda }\left( t\right) =-\sigma _{1}\left( t\right) .
\label{(17)}
\end{equation}%
Finally, choose $\rho \left( t\right) =\widehat{\rho }\left( t\right) $, such
that 
\begin{align}
\widehat{\rho }( t) &=-\Big[ \alpha _{1}( t) +\frac{1}{%
2}\sigma _{1}^{2}( t) -\sigma_1^2(t)+\int_{ \mathbb{R^*} }\Big( e^{\widehat{\theta 
}( t,\zeta) }-1-\widehat{\theta }( t,\zeta) +( e^{%
\widehat{\theta }( t,\zeta) }-1) \gamma _{1}( t,\zeta)
\Big) \nu  d\zeta) \Big]\nonumber\\
&=-\left[ \alpha _{1}\left( t\right) -\frac{1}{%
2}\sigma _{1}^{2}\left( t\right) +\int_{ \mathbb{R^*} }\Big( \ln(1+\gamma_1(t,\zeta)) - \gamma_1(t,\zeta)
\Big) \nu ( d\zeta) \right].\label{(18)}
\end{align}%
Then 
\begin{align}
  \widehat{A}_t:&=A_t^{(\widehat{\rho},\widehat{\lambda},\widehat{\theta})}\nonumber\\
  &=-\Big[\int_0^t\{ \alpha_1(s)-\frac{1}{2}  \sigma_1^2(s) + \int_{ \mathbb{R^*} }(\ln(1+\gamma_1(s,\zeta))-\gamma_1(s,\zeta))\nu(d\zeta)\} ds \nonumber\\
  &+\int_0^t \sigma_1(s)dB(s)+\int_0^t \int_{ \mathbb{R^*} } \ln(1+\gamma_1(s,\zeta))\widetilde{N}(ds,d\zeta) \Big],\label{2.20a}
\end{align}
with $\widehat{Y}_{t}=Y_{t}^{\left( \widehat{\rho },\widehat{\lambda },%
\widehat{\theta }\right)} = \exp(\widehat{A}_t)$ we have, by (\ref{(14)})%
\begin{eqnarray}
d\left( \widehat{Y}_{t}X\left( t\right) \right) -\widehat{Y}_{t}dX\left(
t\right) +\widehat{Y}_{t}X\left( t\right) d\Gamma _{t}  \label{(19)}
=\widehat{Y}_{t}C\left( t\right) dK_{t}.  
\end{eqnarray}%
Substituting this into (\ref{(10)}), we get 
\[
d\left( \widehat{Y}_{t}X\left( t\right) \right) -\widehat{Y}_{t}C\left(
t\right) dK_{t}=\widehat{Y}_{t}C\left( t\right) d\Lambda _{t},
\]%
which we integrate to, since $\widehat{Y}_0=1$, 
\[
\widehat{Y}_{t}X\left( t\right) =z+\int_{0}^{t}\widehat{Y}_{s}C\left(
s\right) d\left( K_{s}+\Lambda _{s}\right).
\]%
Solving for $X\left( t\right) $, we obtain the following:
\begin{theorem}\label{th2.3}
With initial value $z$ the corresponding optimal wealth process $\widehat{X}_z(t)$ is given by
\begin{align}
\widehat{X}_z(t) &=z\widehat{Y}_{t}^{-1}+\widehat{Y}_{t}^{-1}\int_{0}^{t}
\widehat{Y}_{s}C( s) d( K_{s}+\Lambda
_{s})    \nonumber\\
&=z\exp ( -\widehat{A}_t) +\exp ( -\widehat{A}_t)\int_{0}^{t}\exp (\widehat{A}_s)
C\left( s\right) d\left( K_{s}+\Lambda _{s}\right).   \label{(20)}
\end{align}

\end{theorem}
In particular, note that%
\begin{equation}
\frac{d}{dz}\widehat{X}_{z}\left( t\right) =\exp ( -\widehat{A}_{t}) . \label{(21)}
\end{equation}
\subsection{The optimal initial wealth and the option price $\widehat{z}$}
Completing the Stackelberg game, we now proceed to find the initial wealth $\widehat{z}$ which leads to a response $\widehat{\pi}=\pi_{\widehat{z}}$ which is optimal for Problem \ref{mv}, in the sense that  $J(\widehat{z},\widehat{\pi}) \leq J(z,\pi)$ over all pairs $(z,\pi)$.\\
To this end, choose $z \in \mathbb{R}$ and let $\widehat{\pi}_z$ be the corresponding optimal portfolio given by \eqref{3.9a} and let $\widehat{X}_z$ be the corresponding optimal wealth process given by \eqref{2.10} and \eqref{(20)}.  \\
Then
\begin{align*}
    \inf_{z,\pi}J(z,\pi)=\inf_{z,\pi}E\Big[\frac{1}{2}(X_{z,\pi}(T)-F)^2\Big]=\inf_{z}E\Big[\frac{1}{2}(X_{z,\widehat{\pi}_z}(T) -F)^2\Big]=\inf_{z}E\Big[\frac{1}{2}(\widehat{X}_{z}(T) -F)^2\Big].
\end{align*}

Note that, if we define
\begin{align*}
    R_t&:=\exp(-\widehat{A}_t)=\exp\Big[\int_0^t\Big\{ \alpha_1(s)-\frac{1}{2}  \sigma_1^2(s) + \int_{ \mathbb{R^*} }\Big(\ln(1+\gamma_1(s,\zeta))-\gamma_1(s,\zeta)\Big)\nu(d\zeta)\Big\} ds \nonumber\\
  &+\int_0^t \sigma_1(s)dB(s)+\int_0^t \int_{ \mathbb{R^*} } \ln(1+\gamma_1(s,\zeta))\widetilde{N}(ds,d\zeta) \Big],
\end{align*}
and
\small
\begin{align*}
    Z_t^*&:=\exp\Big(-\int_0^t \alpha_1(s)ds\Big)R_t=\exp\Big[\int_0^t\Big\{-\frac{1}{2}  \sigma_1^2(s) + \int_{\mathbb{R^*}}(\ln(1+\gamma_1(s,\zeta))-\gamma_1(s,\zeta))\nu(d\zeta)\Big\} ds \nonumber\\
  &+\int_0^t \sigma_1(s)dB(s)+\int_0^t \int_{ \mathbb{R^*} } \ln(1+\gamma_1(s,\zeta))\widetilde{N}(ds,d\zeta) \Big],
\end{align*}
then we can verify by the It\^{o} formula that
\begin{align*}
    dR_t&=R_t\Big(\alpha_1(t)dt + \sigma_1(t) dB(t) + \int_{ \mathbb{R^*} } \gamma_1(t,\zeta) \widetilde{N}(dt,d\zeta)\Big)\nonumber\\
    &=R_t G(t)\Big(\alpha(t)dt + \sigma(t) dB(t) + \int_{ \mathbb{R^*} } \gamma(t,\zeta) \widetilde{N}(dt,d\zeta)\Big)\nonumber\\
    &= R_t G(t) S^{-1}(t)dS(t),
    \end{align*}
    and
 \begin{align} \label{Z*} 
 dZ^{*}_t&=Z^{*}_t G(t)\Big( \sigma(t) dB(t) + \int_{ \mathbb{R^*} } \gamma(t,\zeta) \widetilde{N}(dt,d\zeta)\Big).
\end{align}
\begin{proposition}
Assume that $Z_t^{*}$ is a $P$-martingale. (See \eqref{KS}.) Define
\begin{align}\label{Q*}
    dQ^{*}(\omega)=Z^{*}_T(\omega) dP(\omega) \text { on } \mathcal{F}_T.
\end{align}
Then $Q^{*}$ is an EMM for $S(\cdot)$.
\end{proposition}\vskip 0.3cm 
\dproof To see this, we verify that the coefficients 
$\theta_0(t):=G(t)\sigma(t)$ and $\theta_1(t,\zeta):=G(t)\gamma(t,\zeta)$ satisfy condition \eqref{emm}:

 \begin{align}\label{emm2}
&\alpha(t)+\theta_{0}(t)\sigma(t)+\int_{
\mathbb{R^*}
}\theta_{1}(t,\zeta)\gamma(t,\zeta)\nu(d\zeta)\nonumber\\
&=\alpha(t)-\frac{\alpha(t)}{\sigma^2(t) + \int_{ \mathbb{R^*} } \gamma^2(t,\zeta) \nu(d\zeta)}\sigma^2(t) -\frac{\alpha(t)}{\sigma^2(t) + \int_{ \mathbb{R^*} }\gamma^2(t,\zeta) \nu(d\zeta)}\int_{ \mathbb{R^*} }\gamma^2(t,\zeta) \nu(d\zeta)\nonumber\\
&= \alpha(t) -\alpha(t)=0\text{, }\quad t\geq 0.
\end{align}   
\fproof

Using this we obtain the following, which is the main result in this section:
\begin{theorem}\label{th2.5}
(i) The unique minimal variance price $\widehat{z}= p_{mv}(F)$ of a European option with terminal payoff $F$ at time $T$ is given by
\begin{align} \label{3.14a}
\widehat{z}= \frac{E\Big[ exp(-\widehat{A}_T)\Big(F -exp(-\widehat{A}_T) \int_0^T exp(\widehat{A}_s)C\left( s\right) d\left( K_{s}+\Lambda _{s}\right)\Big)\Big]}{E[exp(-2\widehat{A}_T)]},
\end{align}
where $\widehat{A}_T$ is given by \eqref{2.20a}, $C(s), \Lambda_s$ are given by \eqref{(2)}, \eqref{(3)} respectively, and $K$ is given by \eqref{(14)}.

(ii) Assume that the coefficients $\alpha(t), \sigma(t)$ and $\gamma(t,\zeta)$ are bounded and deterministic.
Then
\begin{align}
 \widehat{z}=p_{mv}(F)=E_{Q^{*}}[F],  \label{2.43} 
\end{align}
where $Q^{*}$ is the EMM measure given by \eqref{Q*}.
\end{theorem}

\dproof 
(i) To minimize $J_0(z):= E\Big[\frac{1}{2} (\widehat{X}_{z}(T)-F)^2\Big]$ with respect to $z$ we note by \eqref{(20)} and \eqref{(21)} that
\begin{align}
\frac{d}{dz} J_0(z)
&=E\Big[(\widehat{X}_{z}(T) -F)\frac{d}{dz}\widehat{X}_{z}(T)\Big] \nonumber\\
&=E[(\widehat{X}_z(T) -F)exp(-\widehat{A}_T)] \label{2.43a}\\
&=E\Big[\Big(z exp(-\widehat{A}_T) +exp(-\widehat{A}_T)\int_0^T exp(\widehat{A}_s)C\left( s\right) d\left( K_{s}+\Lambda _{s}\right)-F\Big)exp(-\widehat{A}_T)\Big].\nonumber
\end{align}
This is $0$ if and only if \eqref{3.14a} holds.
\vskip 0.2cm
\noindent (ii) 
By \eqref{2.43a}, we get
\begin{align*}
 &E[\widehat{X}_z(T)\exp(-\widehat{A}_T)]=E[F \exp(-\widehat{A}_T)],\nonumber\\
 \text{ i.e. }\nonumber\\
 &E\Big[\widehat{X}_z(T)\exp\Big(\int_0^T \alpha_1(s)ds\Big) Z_T^{*}\Big]=E\Big[F\exp\Big(\int_0^T \alpha_1(s)ds\Big) Z_T^{*}\Big].
\end{align*}
If $\alpha_1$ is deterministic, we can cancel out the factor $\exp\Big(\int_0^T \alpha_1(s)ds\Big)$ and \eqref{2.43} follows.\fproof

\begin{remark}
An important question is:
Is $\widehat{z}$ an arbitrage-free price of $F$?\\
If the coefficients are deterministic, we know that the answer is yes, by Theorem \ref{th2.5} (ii). But in the general case in Theorem \ref{th2.5} (i) this is not clear.
In the next section we will give, under some conditions, an affirmative answer to this question, by proving that 
$$\widehat{z}=E_{Q_0}[F],$$
for any $Q_0 \in \mathbb{M}_0$, where $\mathbb{M}_0$ is a nonempty subset of $\mathbb{M}$. As pointed out in the Introduction this implies in particular that $\widehat{z}$ is arbitrage-free.
\end{remark}

\begin{example} \textbf{(European call option)}
We give some details about how to proceed if we want to compute the minimal variance price $\widehat{z}=p_{mv}(F)$ explicitly in the case of a European call option:
\begin{myenumerate}
    \item 
Note that the term $C(s)$ in Theorem \ref{th2.5} depends on the coefficients $\beta$ and $\kappa$ in the It\^o representation of $F$. These coefficients can for example be found by using the generalised Clark-Ocone formula for L\' evy processes, extended to $L^2(P)$. See Theorem 12.26 in Di Nunno et al (2009). \\
Let us find these coefficients in the case of a European call option, where
$$F=(S(T)-K)^{+},$$
where $K$ is a given exercise price. In this case $F(\omega)$ represents the payoff at time $T$ (fixed) of a (European call) option which gives the owner the right to buy the stock with value $S(T,\omega)$ at a fixed exercise price K. Thus if $S(T,\omega)>K$ the owner of the option gets the profit $S(T,\omega)-K$ and if $S(T,\omega)\leq K$ the owner does not exercise the option and the profit is 0. Hence in this case
$$F(\omega)=(S(T,\omega)-K)^{+}.$$
Thus, we may write
$$F(\omega)=f(S(T,\omega)),$$
where

$$f(x)=(x-K)^{+}.$$

The function $f$ is not differentiable at $x=K$, so we cannot use the chain rule directly to evaluate $D_tF$. However, we can approximate f by $C^1$ functions $f_n$ with the property that 
$$f_n(x)=f(x) \quad \text{ for } \quad |x-K| \geq \frac{1}{n},$$
and 
$$0\leq f'_n(x)\leq 1 \text{ for all } x.$$
Putting
$$F_n(\omega)=f_n(S(T,\omega)),$$
we see
$$D_tF(\omega)=\lim_{n \to +\infty} D_tF_n(\omega).$$

If the coefficients $\alpha, \sigma,\gamma$ of the risky asset price $S$ are deterministic, we get
\begin{align*}
    \beta(t)= E[D_t F | \mathcal{F}_t], 
    \kappa(t,\zeta)=E[D_{t,\zeta}F | \mathcal{F}_t],
\end{align*}
where $D_t F$ and $D_{t,\zeta} F$ denote the generalised Malliavin derivatives (also called the Hida-Malliavin derivative) of F at  $t$ and $(t,\zeta)$ respectively, with respect to $B(\cdot)$ and $N(\cdot,\cdot)$, respectively. Combining this with the chain rule for the Hida-Malliavin derivative and the Markov property of the process $S(\cdot)$, and assuming for simplicity that $\sigma$ is constant and $\gamma(t,\zeta)=\gamma(\zeta)$ does not depend on $t$,  we obtain the following for $\beta$:
\begin{align}\label{beta}
    \beta(t)&= E^{S_0}\Big[\chi_{[K,\infty)}(S(T))\sigma S(T) |\mathcal{F}_t\Big]\nonumber\\
    &=E^{S(t)}\Big[\chi_{[K.\infty)}(S(T-t))\sigma S(T-t)\Big].
\end{align}
To find the corresponding result for $\kappa$ we first use the chain rule for $D_{t,\zeta}$ (Theorem 12.8 in Di Nunno et al (2009)) 
and get
\begin{align*}
    D_{t,\zeta}S(T)&=D_{t,\zeta}\Big[S_0 \exp \Big(\alpha T -\frac{1}{2}\sigma^2 T+ \sigma B(T) +\int_{\mathbb{R^*}}(\log(1+\gamma(\zeta))-\gamma(\zeta))\nu(d\zeta)T\nonumber\\&+\int_0^t\int_{\mathbb{R^*}}\ln(1+\gamma(\zeta))\widetilde{N}(ds,d\zeta) \Big)\Big]=S(T)\gamma(\zeta).
\end{align*}
Then we obtain
\begin{align} \label{kappa}
    \kappa(t,\zeta)&=E^{S_0}\Big[\chi_{[K,\infty)}(S(T)+D_{t,\zeta}S(T))-\chi_{[0,T]}(S(T))\Big|\mathcal{F}_t\Big]\nonumber\\
    &=E^{S_0}\Big[\chi_{[K,\infty)}(S(T)+\gamma(\zeta)S(T))
    -\chi_{[0,T]}(S(T))\Big|\mathcal{F}_t\Big]\nonumber\\&=E^{S(t)}\Big[\chi_{[K,\infty)}(S(T-t)+\gamma(\zeta)S(T-t)) - \chi_{[K,\infty)}(S(T-t))\Big],
\end{align}
where in general $E^{y}[h(S(u))]$ means $E[h(S^{y}(u))]$, i.e.   expectation when $S$ starts at $y$.
\item
Assume that the coefficients $\alpha, \sigma$ and $\gamma$ of the process $S$ are deterministic and bounded. To compute numerically the minimal variance price $$\widehat{z}=p_{mv}((S(T)-K)^{+})=E\Big[(S(T)-K)^{+}Z_T^{*}\Big]$$
of a European call option with payoff  
\begin{align*}
    F=(S(T)-K)^{+}= E[F]+\int_0^T \beta(t)dB(t)+\int_0^T \int_{\mathbb{R^*}} \kappa(t,\zeta) \widetilde{N}(dt,d\zeta),
\end{align*}
where $\beta,\kappa$ are given by \eqref{beta},\eqref{kappa}, respectively, we use the It\^o formula combined with \eqref{Z*} to obtain
\begin{align*}
    \widehat{z}=E[F Z_T^{*}]=E[F]+ \int_0^T G(t)\Big\{\sigma(t) E[Z_t^{*}\beta(t)] + \int_{\mathbb{R^*}}\gamma(t,\zeta) E[\kappa(t,\zeta) Z_t^{*}] \nu(d\zeta)\Big\} dt,
\end{align*}
where $G(t)$ is given by \eqref{(6)}, i.e.
\begin{equation}
G\left( t\right) =-\alpha \left( t\right) \left( \sigma ^{2}\left( t\right)
+\int_{ \mathbb{R^*} }\gamma ^{2}\left( t,\zeta\right) \nu \left( d\zeta\right) \right)
^{-1}.  \label{2.53}
\end{equation}

\item
Alternatively, in some cases it may be convenient to use the Fourier transform in the computation, as follows:\\
Recall that if $\eta_t \in L^2(P)$ is a L\' evy process with the representation
\begin{align*}
    \eta_t=\alpha_0 t + \sigma_0 B(t)+\int_{\mathbb{R^*}}\zeta \widetilde{N}(t,d\zeta),
\end{align*}
where $\alpha_0$ and $\sigma_0$ are constants, then
\begin{align*}
    E[e^{iu\eta_t}]= e^{t\Psi(u)},
\end{align*}
where
\begin{align*}
\Psi(u)=i\alpha_0 u -\frac{1}{2}\sigma_0^2 u^2 +\int_{\mathbb{R^*}}(e^{iu\zeta} - 1 -iu\zeta)\nu(d\zeta),
\end{align*}
$\nu$ being the L\' evy measure of $\eta$.
Combining this with the Fourier transform inversion
\begin{align*}
    f(x)=\frac{1}{\sqrt{2\pi}}\int_{\mathbb{R}}\widehat{f}(y)e^{iyx}dy,
\end{align*}
where in general
\begin{align*}
    \widehat{f}(y)=\frac{1}{\sqrt{2\pi}}\int_{\mathbb{R}}f(x)e^{-iyx}dx
\end{align*}
is the Fourier transform of $f$, we see that we can obtain explicit expressions of the type
\begin{align*}
E[f(\eta_t)]&=E\Big[\frac{1}{\sqrt{2\pi}}\int_{\mathbb{R}}\widehat{f}(y)e^{iy\eta_t}dy\Big]
=\frac{1}{\sqrt{2\pi}}\int_{\mathbb{R}}\widehat{f}(y)E[e^{iy\eta_t}]dy\nonumber\\
&=\frac{1}{\sqrt{2\pi}}\int_{\mathbb{R}}\widehat{f}(y)e^{t\Psi(y)}dy,
\end{align*}
and similarly (by extending to 2 dimensions) for $\widehat{z}=E\Big[(S(T)-K)^{+} Z_T^{*}\Big]$. We omit the details.
\end{myenumerate}
\end{example}

\section{Finding the minimal variance option price by means of the maximum principle and EMMs}

In this section we relate the minimal variance option price $\widehat{z}$  found in Section 2 to pricing by means of EMMs. We will do this by approaching Problem \ref{cost} by means of stochastic control theory. Since the given $\mathcal{F}_T$-measurable random variable $F$ is not Markovian, we cannot use classical dynamic programming to solve it. However, we can adapt the maximum principle method in Agram et al (2020) 
 to our situation. Thus we define the Hamiltonian $H$ by%
\begin{align}
H(t,x,z,\pi,p,q,r)=x\pi\alpha p+x\pi\sigma q+\int_{%
\mathbb{R^*}
}x\pi\gamma(\zeta)r(\zeta)\nu(d\zeta),
\end{align}
and we define the associated adjoint BSDE for the adjoint processes $(p,q,r)$ by%
\begin{align} \label{adj-p}
\left\{
\begin{array}
[c]{ll}%
dp(t) & =-\pi(t)\left[  \alpha(t)p(t)+\sigma(t)q(t)+\int_{%
\mathbb{R^*}
}\gamma(t,\zeta)r(t,\zeta)\nu(d\zeta)\right]  dt\\
& +q(t)dB(t)+\int_{%
\mathbb{R^*}
}r(t,\zeta)\widetilde{N}(dt,d\zeta),\\
p(T) & =X(T)-F.
\end{array}
\right.
\end{align}
Then, by an easy extension of Theorem 2.4 in Agram et al (2020) 
to jumps, we have:
\begin{theorem}
Suppose $\widehat{z},\widehat{\pi}$ is a solution of the problem \eqref{cost}, with associated solutions $(\widehat{p},\widehat{q},\widehat{r})$ of the corresponding BSDE \eqref{adj-p}. Then 
\begin{align*}
&\widehat{p}(0)+\nabla_{z}\widehat{H}(t)=0\text{ at }z=\widehat{z},\\
&\text{ and }\nonumber\\
&\nabla_{\pi}\widehat{H}(t)=0\text{ at }\pi=\widehat{\pi},
\end{align*}
i.e.,%
\begin{align}
&\widehat{p}(0)=0,\label{p_0}\\
&\text{ and }\nonumber\\
&\widehat{X}(t)\left[  \alpha(t)\widehat{p}(t)+\sigma(t)\widehat{q}(t)+\int_{%
\mathbb{R^*} \label{cond}
}\gamma(t,\zeta)\widehat{r}(t,\zeta)\nu(d\zeta)\right]  =0.
\end{align}
\end{theorem}
Note that if $z=0$ then $X(t)=0$ for all $\pi$ and all $t$, which is a trivial special case. Thus we assume that $z \neq 0$ and this gives $X(t) \neq 0$ for all $\pi$ and $t$ and then \eqref{p_0} and \eqref{cond} give 
\begin{align} \label{cond1}
 \alpha(t)\widehat{p}(t)+\sigma(t)\widehat{q}(t)+\int_{
\mathbb{R^*}
}\gamma(t,\zeta)\widehat{r}(t,\zeta)\nu(d\zeta) =0.
\end{align}
and 
\begin{align}  \label{adjp2}
d\widehat{p}(t) & =\widehat{q}(t)dB(t)+\int_{
\mathbb{R^*}
}\widehat{r}(t,\zeta)\widetilde{N}(dt,d\zeta),\quad \widehat{p}(0)=0 \text{ and }  
\widehat{p}(T) =\widehat{X}(T)-F.
\end{align}
We have proved:
\begin{theorem} \label{4.2}
Suppose there exists an optimal control $(\widehat{z},\widehat{\pi})$ of the problem \eqref{cost}.
Then the corresponding forward-backward system, consisting of 
\begin{align}
d\widehat{X}(t) & =\widehat{X}(t)\widehat{\pi}(t)\left[  \alpha(t)dt+\sigma(t)dB(t)+\int_{
\mathbb{R^*} \label{system}
}\gamma(t,\zeta)\widetilde{N}(dt,d\zeta)\right], \quad \widehat{X}(0)=\widehat{z},\\
d\widehat{p}(t) & =\widehat{q}(t)dB(t)+\int_{
\mathbb{R^*}}\widehat{r}(t,\zeta)\widetilde{N}(dt,d\zeta),\quad \widehat{p}(T) =\widehat{X}(T)-F,
\label{adjoint}
\end{align}
satisfies the equations
\begin{align*}
    \widehat{p}(0)&=0,\\
    \alpha(t)\widehat{p}(t)+\sigma(t)\widehat{q}(t)+\int_{\mathbb{R^*}}\gamma(t,\zeta)\widehat{r}(t,\zeta)\nu(d\zeta) &=0.
\end{align*}
\end{theorem}
 This implies the following option pricing result:
 \begin{theorem}{\bf{(Option pricing theorem 1)}}\label{th2.2}
 Suppose $(\widehat{z},\widehat{\pi})$ is an optimal control for the problem \eqref{cost}, with corresponding solutions $\widehat{X}=X_{\widehat{z},\widehat{\pi}}$ and $(\widehat{p},\widehat{q},\widehat{r})$ of
 \eqref{system} and \eqref{adjoint}, respectively.
 Then for all $Q \in \mathbb{M}$ the minimal variance price of $F$ is given by
 \begin{align} \label{pricing1}
     \widehat{z}=p_{mv}(F)= E_Q[F] +E_Q \Big[\int_0^T \widehat{q}(t)dB(t) + \int_0^T \int_{\mathbb{R^*}} \widehat{r}(t,\zeta)\widetilde{N}(dt,d\zeta)\Big], \text{ for \bf{all} } Q \in \mathbb{M}.
 \end{align}
 
 \end{theorem}
 \dproof
 By \eqref{adjoint} in Theorem  \ref{4.2} we have, for all $Q \in \mathbb{M}, $
 \begin{align*}
 E_Q[\widehat{p}(T)]
 =E_Q[\widehat{X}(T)- F] =E_Q[\widehat{X}(T)]-E_Q[F]=\widehat{z}-E_Q[F].   
 \end{align*}
 Hence, by \eqref{adjoint},
 \begin{equation*}
     \widehat{z}=E_Q[F]+E_Q[\widehat{p}(T)]=
     E_Q[F] +E_Q\Big [\int_0^T \widehat{q}(t)dB(t) + \int_0^T \int_{\mathbb{R^*}} \widehat{r}(t,\zeta)\widetilde{N}(dt,d\zeta)\Big].
 \end{equation*}
 \fproof\\
 
 From this we deduce the following result, which shows in particular that if the market is complete, then the minimal variance price agrees with the price given by the classical Black-Scholes formula:
 \begin{corollary}{\bf{(Generalised Black-Scholes formula)}}\label{cor2.3}
 Let $\widehat{X}$ be as in Theorem \ref{4.2}. Suppose that
 \begin{align*}
 \widehat{X}(T)=F \text{ a.s. }
 \end{align*}
 Then
 \begin{align*}
     \widehat{z}=p_{mv}(F)= E_Q[F], \text{ for \emph{all} } Q \in \mathbb{M}.
 \end{align*}
 \end{corollary}
 \dproof\\
 If $\widehat{X}(T)=F \text{ a.s. }$ then clearly $(\widehat{z},\widehat{\pi})$ is an optimal pair for the problem \eqref{cost}. Therefore, by \eqref{adjoint} we have $\widehat{p}(T)=0$ and hence, since $\widehat{p}$ is an $(P,\mathbb{F})$-martingale,
 \begin{align*}
     p(t)=E[p(T) | \mathcal{F}_t]=0,
 \end{align*}
 for all $t$. But then $\widehat{q}(t)=\widehat{r}(t,\zeta) = 0$ for all $t,\zeta$, and the result follows from \eqref{pricing1}.
 \fproof
 
We now make the following definition:
\begin{definition}
    $\mathbb{M}_0$ is the set of measures $Q \in \mathbb{M}$ which are also equivalent martingale measures for the process $\widehat{p}(\cdot)$. 
\end{definition}

Then we obtain the following pricing result:
\begin{theorem}{\bf{(Option pricing theorem 2)}}
    Let $\widehat{z},\widehat{\pi}, \widehat{p},\widehat{q},\widehat{r}$ be as in Theorem \ref{th2.2}. Then we have 
    \begin{align*}
        \widehat{z}=p_{mv}(F)= E_{Q_0}[F],
    \end{align*}
    for all $Q_0 \in \mathbb{M}_0$.
\end{theorem}
\dproof
If we apply \eqref{pricing1} to $Q_0 \in \mathbb{M}_0$ we get $\widehat{z}=p_{mv}(F)=E_{Q_0}[F]$, because\\
$$E_{Q_0}\Big[\int_0^T \widehat{q}(t)dB(t) + \int_0^T \int_{\mathbb{R^*}} \widehat{r}(t,\zeta)\widetilde{N}(dt,d\zeta)\Big]=0.$$
\fproof\\
The following result illustrates how the measures in $\mathbb{M}_0$ may look like: 
\begin{theorem}{\bf{(Option pricing theorem 3)}}\label{th3.7}
Suppose the risky asset price $S(t)$ is given as in \eqref{S}, but with $m=k=1$, i.e. 
\begin{align}\label{S0}
\begin{array}
[c]{ll}%
dS(t) & =S(t)\left[  \alpha(t)dt+\sigma(t)dB(t)+\int_{%
\mathbb{\mathbb{R^*}}
}\gamma(t,\zeta)\widetilde{N}(dt,d\zeta)\right],\quad S(0) >0.
\end{array}
\end{align}
Specifically, suppose that $(\widehat{\theta}_0(t),\widehat{\theta}_1(t,\zeta))$ solves the following system of two equations
\begin{itemize}
    \item
    \begin{align}
  \alpha(t)+\widehat{\theta}_{0}(t)\sigma(t)+\int_{\mathbb{R^*}
}\widehat{\theta}_{1}(t,\zeta)\gamma(t,\zeta)\nu(d\zeta)=0\text{, } \quad t\geq 0,\label{2.21a}
\end{align}
\item
\begin{align}
\widehat{\theta}_{0}(t)\widehat{q}(t)+\int_{\mathbb{R^*}
}\widehat{\theta}_{1}(t,\zeta)\widehat{r}(t,\zeta)\nu(d\zeta)=0\text{, }\quad t\geq 0.\label{2.22a}
\end{align}
\end{itemize}
Define the process $\widehat{Z}(t)=Z^{\widehat{\theta}_{0},\widehat{\theta}_{1}}(t)$ by
\begin{align*}
d\widehat{Z}(t)  =\widehat{Z}(t)\Big[\widehat{\theta}_{0}(t)dB(t)+\int_{\mathbb{R^*}
}\widehat{\theta}_{1}(t,\zeta)\widetilde{N}(dt,d\zeta)\Big],\quad 
\widehat{Z}(0) =1,
\end{align*}
i.e.,
\begin{align*}
\widehat{Z}(t)  &  =\exp\left(  \int_{0}^{t}\widehat{\theta}_{0}(s)dB(s)-\frac{1}{2}\int_{0}%
^{t}\widehat{\theta}_{0}^{2}(s)ds\right. \\
&  +\int_{0}^{t}\int_{%
\mathbb{R^*}
}\left\{  \ln(1+\widehat{\theta}_{1}(s,\zeta))-\widehat{\theta}_{1}(s,\zeta)\right\}  \nu
(d\zeta)ds\left.  +\int_{0}^{t}\int_{%
\mathbb{R^*}
}\ln(1+\widehat{\theta}_{1}(s,\zeta))\widetilde{N}(ds,d\zeta)\right).
\end{align*}
Assume that $\widehat{Z}$ is a martingale. (See Kallsen \& Shiryaev (2002)).
Then the measure $\widehat{Q}_0:=Q^{\widehat{\theta}_{0},\widehat{\theta}_{1}}$ defined by%
\begin{align*}
d\widehat{Q}_0(\omega)=\widehat{Z}(T)dP(\omega)\text{ on }\mathcal{F}_{T}
\end{align*}
is in $\mathbb{M}_0$, and
$$p_{mv}(F) =\widehat{z}=  E_{\widehat{Q}_0}[F].$$
\end{theorem}

\section{Examples}
In this section we illustrate the results above by considering some examples.
\subsection{The classical Black-Scholes market}
We consider the classical Black-Scholes market,  with $N=0,m=1$ in the market model \eqref{S}, \eqref{W}. Then  we have a complete financial market with the following two investment possibilities:
\begin{description}
\item[(i)] A risk free asset, with unit price $S_{0}(t)=1$ for all $t$.

\item[(ii)] A risky asset, with unit price $S(t)$ at time $t$ given by%
\[
\begin{array}
[c]{ll}%
dS(t)  =S(t)\left[  \alpha(t)dt+\sigma(t)dB(t)\right],\quad S(0) >0.
\end{array}
\]

\end{description}
Note that with $N(t,\zeta)=0$ and $\sigma(t)$ bounded away from 0 for all $t$,  \eqref{emm} gets the form%

$
\alpha(t)+\theta_{0}(t)\sigma(t)=0\text{,}%
$
which has the unique solution%
$
\theta_{0}(t)=-\frac{\alpha(t)}{\sigma(t)}\text{.}%
$
This gives, by \eqref{Z2},
$
Z(t)=Z^{\theta_{0}}(t)=\exp\left(  \int_{0}^{t}\theta_{0}(s)dB(s)-\frac{1}{2}\int_{0}^{t}\theta_{0}^2(s)ds\right)
$
and, assuming that $Z$ is a martingale,
$
dQ^{\theta_{0}}(\omega)=Z^{\theta_{0}}(T)dP(\omega)
$
is the only element in $\mathbb{M}$. Therefore the unique option price in this case is
\begin{equation*}
 p_b(F)=p_s(F)= E_{Q^{\theta_{0}}}[F],
\end{equation*}
which is the celebrated Black-Scholes formula.\\

We now compare this with what we get by using the stochastic control approach of Section 4:
Since, $N(t,.)=0$ for all $t$, \eqref{2.21a} gets the form%

\begin{equation}
\alpha(t)+\theta_{0}(t)\sigma(t)=0 \text{ , i.e.}
-\frac{\alpha(t)}{\sigma(t)}=\theta_0(t).
\end{equation}

From \eqref{2.22a}, we have 
\[
\widehat{p}(t)=\widehat{q}(t)=0.
\]
Therefore, by \eqref{p_0} and \eqref{adjp2}, we get
\[
\widehat{X}(T)-F=\widehat{p}(T)=0.
\]
By (\ref{p_0})  this gives
\[
\widehat{X}(T)=X^{\widehat{z},\widehat{\pi}}(T)=F.
\]
Since $Q^{\theta_0}$ is a martingale measure for $\widehat{X}(t)$, we conclude that
\[
\widehat{z}=X^{\widehat{z},\widehat{\pi}}(0)=E_{Q^{\theta_0}}[X^{\widehat{z},\widehat{\pi}}(T)]
=E_{Q^{\theta_0}}[F].
\]
We conclude that in this case the optimal control $\widehat{\pi}(t)$ is the replicating portfolio for $F$, and the optimal initial wealth $\widehat{z}\in\mathbb{R}$ is the $Q^{\theta_0}$-expectation of $F$. Thus the minimum variance price $p_{mv}(F)$ coincides with the classical option price in this case.

\subsection{A continuous incomplete market}
Consider the case with no jumps ($N=0$) and with two Brownian motions, $B_1(t), B_2(t)$. Then the price process is given by
\begin{align}
dS(t)=S(t)[\alpha(t)dt +\sigma_1(t)dB_1(t) + \sigma_2(t)dB_2(t)],\quad \label{2.25}
S(0) > 0, 
\end{align}
and equations \eqref{2.21a},\eqref{2.22a} become
   \begin{align}\label{2.27}
  \alpha(t)+x_1 \sigma_1(t)+x_2 \sigma_2(t)&=0\text{, }\quad t\geq0, \\
x_1 \widehat{q}_1(t)+x_2 \widehat{q}_2(t)&=0\text{, }\quad t\geq0.\label{2.28}
    \end{align}
where we for simplicity have put $\widehat{\theta_0}(t)=(x_1,x_2)$. This is a linear system of two equations with the two unknowns $x_1, x_2$. This system has a unique solution $(x_1,x_2)=\widehat{\theta_0}(t)$ if and only if
$\sigma_1(t)\widehat{q}_2(t) - \sigma_2(t) \widehat{q}_1(t) \neq 0.$ We conclude that
\begin{corollary}
(a) In the market \eqref{2.25} there is a unique $\widehat{Q}_{0} \in \mathbb{M}_0$ if and only if 
$\sigma_1(t)\widehat{q}_2(t) - \sigma_2(t) \widehat{q}_1(t) \neq 0$, and if the process $Z=\widehat{Z}_0$ defined by \eqref{Z1} is a martingale, then the minimal variance price of F is given by
\begin{align*}
    p_{mv}(F)=E_{\widehat{Q}_{0}}[F], \quad \text{ with  } d\widehat{Q}_0= Z^{\widehat{\theta}_0}(T)dP,
\end{align*}
where $\widehat{\theta}_0=(x_1,x_2)$ is the unique solution of the system \eqref{2.27}-\eqref{2.28}.

(b) If the coefficients $\alpha,\sigma_1,\sigma_2$ are deterministic, we can apply Theorem \ref{th2.5} (ii) to conclude that
\begin{align*}
    p_{mv}(F)=E_{Q^{*}}[F],
\end{align*}
where $dQ^{*}=Z^{*}(T) dP$ is given by \eqref{Q*}.
\end{corollary}
\subsection{A pure jump incomplete market}
Suppose the risky asset price $S(t)$  is given by
\begin{align}\label{purejump}
    dS(t) = S(t)[\alpha(t) dt + \int_{\mathbb{R^*}} \gamma(t,\zeta)\widetilde{N}(dt,d\zeta)],\quad 
    S(0) > 0.
\end{align}
Then the wealth process satisfies
\begin{align}
d\widehat{X}(t) & =\widehat{X}(t)\widehat{\pi}(t)\left[  \alpha(t)dt+\int_{
\mathbb{R^*} \label{drif-jum}
}\gamma(t,\zeta)\widetilde{N}(dt,d\zeta)\right], \quad \widehat{X}(0)=\widehat{z}.
\end{align}
Here there is no Brownian motion component and only one compensated Poisson random measure (i.e. $k=1$), but we assume that $N$ has at least two possible jump sizes, i.e. that the L\' evy measure $\nu$ is not a point mass. Then the market is not complete, because there are several solutions $\theta_1(t,\zeta)$ of the equation \eqref{emm} (or \eqref{2.21a}), which now has the form
\begin{equation}\label{2.32}
 \alpha(t) + \int_{\mathbb{R^*}} \theta_1(t,\zeta) \gamma(t,\zeta) \nu(d\zeta) = 0.   
\end{equation}
To find possible elements $\widehat{Q}$ of $\mathbb{M}_0$ we combine
\eqref{2.32} with \eqref{2.21a}, which now reduces to
\begin{align} \label{2.33}
    \int_{\mathbb{R^*}} \theta_1(t,\zeta)\widehat{r}(t,d\zeta) \nu(d\zeta)=0.
\end{align}
This gives the following result:
\begin{corollary}
(a) Suppose there exists a solution $\theta_1(t,\zeta)=\widehat{\theta}_1(t,\zeta)$ of the two equations \eqref{2.32}, \eqref{2.33} and that the corresponding $Z^{\widehat{\theta}_1}$ defined by \eqref{Z1} is a martingale. Then the minimal variance price of $F$ in the market \eqref{purejump} is given by
\begin{equation*}
    p_{mv}(F)= E_{\widehat{Q}}[F], \quad \text{ where  }
    d\widehat{Q}=Z^{\widehat{\theta}_1}(T) dP.
\end{equation*}
(b) If the coefficients $\alpha,\gamma$ are deterministic, we can apply Theorem \ref{th2.5} (ii) to conclude that
\begin{align}
    p_{mv}(F)=E_{Q^{*}}[F]
\end{align}
where $dQ^{*}=Z^{*} dP$ is given by \eqref{Q*}.

\end{corollary}

\subsection{Merton type markets}
Finally, consider the Merton type markets, driven by a Brownian motion $B(t)$ and a jump process being the standard Poisson process $\mathbf{N}(t)$ with intensity $\lambda >0$, which implies that the L\' evy measure $\nu$ is just the point mass at 1, $\delta_{1}$. Then the corresponding compensated Poisson random measure will be
\begin{align*}
\tilde{N}(dt,d\zeta)
=\delta_{1}(d\zeta)d\mathbf{N}(t)-\lambda \delta_{1}(d\zeta)dt,
\end{align*}
where $\delta_1(d\zeta)$ is the unit point mass at 1,
and the price process is then given by
\begin{align}\label{Merton}
    dS(t) = S(t)[\alpha(t) dt + \sigma(t) dB(t) +  \gamma(t,1)(d\mathbf{N}(t)-\lambda dt)],\quad 
    S(0) > 0.
\end{align}
with corresponding wealth process
\begin{align}
d\widehat{X}(t) & =\widehat{X}(t)\widehat{\pi}(t)\left[  \alpha(t)dt+\sigma (t)dB(t)+ 
 \gamma(t,1)(d\mathbf{N}(t)-\lambda dt)\right], \quad \widehat{X}(0)=\widehat{z}.
\end{align}
Then the equations \eqref{2.21a}, \eqref{2.22a}  get the form
   \begin{align}\label{2.39}
  \alpha(t)+x(t)\sigma(t)+y(t)\gamma(t,1)\lambda &=0\text{, }t\geq0, \\
x(t)\widehat{q}(t)+
y(t)\widehat{r}(t,1)\lambda &=0\text{, }t\geq0,\label{2.40}
    \end{align}
    where we for simplicity have put $x(t)=\widehat{\theta}_0(t),y(t)=\widehat{\theta}_1(t,1).$
      This system has a unique solution $(x(t),y(t))$ if and only if
    \begin{align*}
    \sigma(t)\widehat{r}(t,1) -\widehat{q}(t) \gamma(t,1) \neq 0.
    \end{align*}
    \begin{corollary}
(a) Suppose there exists a solution $x(t)=\widehat{\theta}_0(t),y(t)=\widehat{\theta}_1(t,\zeta)$ of the two equations \eqref{2.39},\eqref{2.40} and that the corresponding $Z^{\widehat{\theta}_0, \widehat{\theta}_1}$ defined by \eqref{Z1} is a martingale. Then the minimal variance price of $F$ in the market \eqref{Merton} is given by
\begin{equation*}
    p_{mv}(F)= E_{\widehat{Q}}[F], \quad \text{ where  }
    d\widehat{Q}=Z^{\widehat{\theta}_0, \widehat{\theta}_1}(T) dP.
\end{equation*}
(b) If the coefficients $\alpha,\sigma,\gamma$ are deterministic, we can apply Theorem \ref{th2.5} (ii) to conclude that
\begin{align*}
    p_{mv}(F)=E_{Q^{*}}[F],
\end{align*}
where $dQ^{*}=Z^{*}(T) dP$ is given by \eqref{Q*}.

\end{corollary}
  
\subsubsection{Pure jump market}
\begin{align} \label{S1}
    dS(t) = S(t)[\alpha(t) dt +  \gamma(t,1)(d\mathbf{N}(t)-\lambda dt)],\quad
    S(0) > 0,
\end{align}
with corresponding wealth process
\begin{align*}
d\widehat{X}(t) & =\widehat{X}(t)\widehat{\pi}(t)\left[  \alpha(t)dt+ 
 \gamma(t,1)(d\mathbf{N}(t)-\lambda dt)\right], \quad \widehat{X}(0)=\widehat{z}.
\end{align*}
Then the equations \eqref{2.39}, \eqref{2.40}  get the form
\begin{align*}
  \alpha(t)+y(t)\gamma(t,1)\lambda &=0 \text{, }\quad t\geq0, \\
y(t)\widehat{r}(t,1)\lambda &=0 \text{, }\quad t\geq0.
\end{align*}
    This system has a unique solutions $\widehat{r}(t,1)=0,y(t)=\widehat{\theta}_1(t,1)=- \frac{ \alpha(t)}{\gamma(t,1)\lambda}$.
    Assume that
    \begin{align*}
        \frac{\alpha(t)}{\gamma(t,1)\lambda} < 1,\quad \text{ for all } t.
    \end{align*}
Defne
\begin{align*}
    dQ(\omega)=Z(T,\omega)dP(\omega) \text{ on } \mathcal{F}_T,
\end{align*}
where
\begin{align*}
Z(t)&= \exp\Big(\int_0^t \lambda\{\ln(1+\widehat{\theta}_1(s,1))-\widehat{\theta}_1(s,1)\}ds+\int_0^t\ln(1+\widehat{\theta}_1(s,1))(d\mathbf{N}(s)-\lambda ds) \Big)\\
&=\exp\Big(\int_0^t -\lambda \widehat{\theta}_1(s,1)ds+\int_0^t\ln(1+\widehat{\theta}_1(s,1))d\mathbf{N}(s) \Big) .
\end{align*}
Suppose that $Z$ is a martingale. Then $Q$ is the unique EMM for the process \eqref{S1}. Hence the unique no-arbitrage price of an option with payoff $F$ in this market is
\begin{align*}
p_{mv}(F)=E_{Q}[F].
\end{align*}
\subsubsection{Merton mixed type market}
Assume that the risky asset price $S(t)$ is given by
\begin{align*}
    dS(t)=S(t)\Big[\alpha_0 dt + \sigma_0 dB(t) + \int_{\mathbb{R^*}} \gamma_0 \widetilde{N}(dt,d\zeta)\Big], \quad S(0)=S_0 > 0,
\end{align*}
where $\widetilde{N}(dt,d\zeta)=N(dt,d\zeta) -\nu(d\zeta)dt$, 
with $N(dt,d\zeta)=d\mathbf{N}(t)\delta_1(d\zeta)$. Here $\mathbf{N}$ is a Poisson process with intensity $\lambda > 0$, $\delta_1$ is the Dirac measure (unit point mass) at 1 and $\alpha_0, \sigma_0 > 0$ and $\gamma_0 > -1$ are given constants. This equation has the following explicit solution (see e.g. Example 1.15 in \O ksendal \& Sulem (2019)):
\begin{align*}
S(t)& =S_0 \exp\Big( \{\alpha_0 -\frac{1}{2}\sigma_0^2-\lambda \gamma_0\}t+\sigma_0 B(t)+\ln(1+\gamma_0)\mathbf{N}(t) \Big).
\end{align*}
In this case the process $Z^{*}_t$ given by \eqref{Z*} gets the form
\begin{align*}
    dZ_t^{*}= Z_{t}^{*}\Big[G\sigma_0 dB(t)+G\int_{\mathbb{R^*}}\gamma_0 \widetilde{N}(dt,d\zeta)\Big], \quad Z_0^{*}=1,
\end{align*}
 $$G\gamma_0>-1,$$
where, by \eqref{(6)},
\begin{align*}
    G= - \frac{\alpha_0}{\sigma_0^2 + \lambda \gamma_0^2}.
\end{align*}
This equation has the solution
\begin{align*}
    Z^{*}_t=\exp \Big(\{-\frac{1}{2}G^2\sigma_0^2 -\lambda G \gamma_0\}t+G \sigma_0 B(t)+ \ln(1+G \gamma_0 )\mathbf{N}(t) \Big).
\end{align*}
By Theorem \ref{th2.5} (ii) the minimal variance price $\widehat{z}$ of a contract with payoff $F$ at time $T$ is
\begin{align*}
    \widehat{z}=E_{Q^{*}}(F)= E[F Z_T^{*}].
\end{align*}
Let us assume that $F$ has the It\^o representation
\begin{align*}
F=F_0+\phi_0 B(T)+\psi_0\mathbf{N}(T), \quad F_0>0,
\end{align*}
where $\phi_0,\psi_0$ are constants.
Then by integration by parts and the It\^o isometry we get that
\begin{align*}
   \widehat{z}= E[F Z_T^{*}]=F_0 + \Big(G \sigma_0 \phi_0 + \lambda [1+G\gamma_0]\psi_0 \Big)T.
\end{align*}

\section{Acknowledgements}
We are grateful to Markus Hess for his valuable comments.

\end{document}